\documentclass{amsart}
\usepackage{amssymb, epsfig}

\newtheorem{thm}{Theorem}[section]

\newtheorem{prop}[thm]{Proposition}

\theoremstyle{definition}
\newtheorem{defn}[thm]{Definition}
\newtheorem{remarkk}[thm]{Remark}
\newenvironment{remark}{\begin{remarkk}}{\end{remarkk}}

\numberwithin{equation}{section}

\renewcommand{\bibname}{\sc}
\renewcommand{\and}{\textnormal{and}\ }
\newenvironment{proof*}{\begin{proof}}{\end{proof}}

\newcommand{\calf}{\mathcal{F}}
\newcommand{\calc}{\mathcal{C}}

\newcommand{\calk}{\mathcal{K}}

\newcommand{\caln}{\mathcal{N}}
\newcommand{\calu}{\mathcal{U}}
\newcommand{\bbz}{\mathbb{Z}}
\newcommand{\bbr}{\mathbb{R}}
\newcommand{\bbq}{\mathbb{Q}}

\newcommand{\bbn}{\mathbb{N}}

\title[Polynomial splittings of metabelian von Neumann
rho--invariants]{Polynomial splittings of metabelian von Neumann
rho--invariants of knots}
\author{Se-Goo Kim}
\address{Department of Mathematics and Research Institute for Basic
Sciences, Kyung Hee University, Seoul 130--701, Korea}
\email{sgkim@khu.ac.kr}
\urladdr{web.khu.ac.kr/\~{}sekim}

\author{Taehee Kim}
\address{Department of Mathematics, Konkuk
University, Seoul 143--701, Korea}
\email{tkim@konkuk.ac.kr} \urladdr{konkuk.ac.kr/\~{}tkim}
\def\subjclassname{\textup{2000} Mathematics Subject Classification}
\expandafter\let\csname subjclassname@1991\endcsname=\subjclassname
\expandafter\let\csname subjclassname@2000\endcsname=\subjclassname
\date{\today}

\begin{document}

\begin{abstract}
    We show that if the connected sum of two knots with coprime Alexander
    polynomials has vanishing von Neumann $\rho$--invariants associated
    with certain metabelian representations then so do both knots.
    As an application, we give a new example of an infinite family of knots
    which are linearly independent in the knot concordance group.
\end{abstract}

\maketitle


\section{Introduction}
A knot $K$ in the 3--sphere $S^3$ is said to be \emph{slice} if
there is a locally flat 2--disk $D$ embedded in the 4--ball $B^4$
with $\partial (B^4,D)=(S^3,K)$. A pair of knots $K_1$ and $K_2$ are
\emph{concordant} if $K_1\# (-K_2)$ is slice where $-K$ is the
mirror image of $K$ with reversed orientation. The set of
concordance classes of knots forms an abelian group under connected
sum, called the \emph{knot concordance group} and denoted by
$\calc$. In $\calc$ the identity 0 is represented by slice knots.

By constructing sliceness obstructions using Seifert forms, Levine
showed that the knot concordance group surjects to $\bbz^\infty
\oplus (\bbz/2)^\infty \oplus (\bbz/4)^\infty$ \cite{lev69a,
lev69b}. A knot with vanishing Levine obstructions is called
\emph{algebraically slice}. Using their own invariants, Casson and
Gordon showed that there are non-slice knots which are algebraically
slice \cite{cg86}. Using Casson--Gordon invariants Jiang showed that
the concordance group of algebraically slice knots is infinitely
generated \cite{jia81}. Gilmer refined Casson--Gordon invariants by
combining Casson--Gordon invariants with the Levine obstructions
\cite{gil93}.

In \cite{cot03} Cochran, Orr and Teichner made progress by
establishing a geometric filtration of $\calc$
\[
0\subset \cdots\subset \calf_{n.5}\subset \calf_{n} \subset \cdots\subset
\calf_{1.5}\subset \calf_{1} \subset \calf_{0.5}\subset
\calf_{0}\subset \calc,
\]
where the subgroup $\calf_{h}$ is the set of all
\emph{$(h)$--solvable} knots. They showed that this filtration is
closely related to the known concordance invariants. For instance, a
knot lies in $\calf_{0.5}$ if and only if the knot is algebraically
slice \cite[Remark 1.3.2]{cot03}. Also they proved that all
previously known concordance invariants vanish on $(1.5)$--solvable
knots \cite[Section 9]{cot03}. In particular, Casson--Gordon--Gilmer
invariants~\cite{cg86, gil93} vanish on $(1.5)$--solvable knots.
Cochran, Orr and Teichner~\cite{cot03, cot04} used von Neumann
$\rho$--invariants ($L^2$--signature defects) to prove that the
quotient group $\calf_2/\calf_{2.5}$ has infinite rank. The second
author \cite{tkim04} proved that $\calf_{1}/\calf_{1.5}$ has an
infinite rank subgroup of knots for which Casson--Gordon invariants
vanish. In their respective papers Cochran--Orr--Teichner
\cite{cot03} and the second author \cite{tkim04} constructed their
examples using genetic modification and showed the linear
independence using signature functions whose integrals over $S^1$
are linearly independent.

Linear independence of knots in $\calc$ may be checked in a
different manner using relative primeness of Alexander polynomials.
Such an approach was first given by Levine~\cite{lev69a} who showed
that if the connected sum of two knots with coprime Alexander
polynomials has vanishing Levine obstructions, then so do both
knots. The first author~\cite{kim05} showed that the
Casson--Gordon--Gilmer invariants split in this way as well.

In this paper, we prove a similar splitting property for the von
Neumann $\rho$--invariants of knots associated with certain
metabelian representations. Let $K$ be a knot and $M_K$ the zero
surgery on the knot $K$ in $S^3$. Let $\Lambda:=\bbq[t^{\pm 1}]$.
For every $x$ in the (rational) Alexander module $H_1(M_K;\Lambda)$,
one can obtain the real-valued $\rho$--invariant $\rho(K,\phi_x)$
where $\phi_x : \pi_1(M_K) \to \bbq(t)/\Lambda \rtimes \bbz$ is the
homomorphism associated with $x$ via the Blanchfield linking form
(see Definition~\ref{defn:representation}). We say that $K$ has
\emph{vanishing $\rho$--invariants} if there exists a
self-annihilating $\Lambda$--submodule $P$ of $H_1(M_K;\Lambda)$
with respect to the Blanchfield linking form (hence algebraically
slice) such that $\rho(K,\phi_x) = 0$ for all $x\in P$. For the
definition of self-annihilating submodule, see
Section~\ref{sec:preliminaries}. Cochran, Orr and Teichner showed
that $(1.5)$--solvable knots have vanishing $\rho$--invariants
\cite[Theorem 4.6]{cot03}. This yields a sliceness obstruction since
a slice knot is $(1.5)$--solvable.

We state the main theorem:

\begin{thm} \label{thm:main}
    Let $K_1$ and $K_2$ be knots with coprime Alexander polynomials. If
    $K_1\# K_2$ has vanishing $\rho$--invariants, then so do
    both $K_1$ and $K_2$.
\end{thm}

\noindent
We give a stronger form of this theorem in Theorem~\ref{thm:splitting}.

To demonstrate the strength of this result,
in Section~\ref{sec:examples}
we give a new example of infinitely many knots with vanishing
Casson--Gordon invariants which are linearly
independent in $\calf_1/\calf_{1.5}$ and hence in $\calc$.

\section{Preliminaries}
\label{sec:preliminaries}

In this section we briefly review the machinery necessary for this
paper. In~\cite{cot03}, Cochran, Orr and Teichner established a
filtration $\{\calf_h\}_{h\in\frac{1}{2}\bbn_0}$ of $\calc$ indexed
by nonnegative half-integers where $\calf_h$ is the subgroup of
\emph{$(h)$--solvable} knots which is defined below. Recall that for
a group $G$ and a nonnegative integer $n$, the \emph{$n$--th derived
group of $G$}, $G^{(n)}$, is defined inductively by the relations
$G^{(0)}:=G$ and $G^{(k)}:=[G^{(k-1)},G^{(k-1)}]$ for $k\ge 1$.

For a CW-complex $W$, we define $W^{(n)}$ to be the regular covering
corresponding to the subgroup $\left(\pi_1(W)\right)^{(n)}$. Suppose
$W$ is an oriented 4--manifold. Then there is an intersection form
$$
\lambda_n : H_2(W^{(n)})\times H_2(W^{(n)}) \rightarrow
\bbz\left[\pi_1(W)/\pi_1(W)^{(n)}\right].
$$
Also there is a self-intersection form $\mu_n$ on $H_2(W^{(n)})$.
For more details about these forms refer to \cite{wal99} and
\cite[Section 7]{cot03}. For a nonnegative integer $n$, an
\emph{$(n)$--Lagrangian} is a submodule $L \subset H_2(W^{(n)})$ on
which $\lambda_n$ and $\mu_n$ vanish and which maps onto a
Lagrangian of $\lambda_0$ under the homomorphism induced by the
covering map.

\begin{defn}\cite[Section 8]{cot03}
    Let $n\in \bbn_0$.
    A knot $K$ is called \emph{$(n)$--solvable} if $M_K$ bounds a spin
    4--manifold $W$ such that the inclusion map $M_K\to W$ induces an
    isomorphism on the first homology and such that $W$ admits an
    $(n)$--Lagrangian with $(n)$--duals. This means that the
    intersection form $\lambda_n$ pairs the $(n)$--Lagrangian
    and the $(n)$--duals nonsingularly and that their images together
    freely generate $H_2(W)$. The 4--manifold $W$ is called an
    \emph{$(n)$--solution for} $K$ and we say $K$ is
    \emph{$(n)$--solvable via} $W$.

    Similarly, we define \emph{$(n.5)$--solvable} knots for $n\in
    \bbn_0$ in such a way that an $(n.5)$--solution $W$ is required to
    admit an $(n+1)$--Lagrangian with $(n)$--duals. For more details,
    refer to \cite[Definitions 8.5 and 8.7]{cot03}.
\end{defn}

Cochran, Orr and Teichner showed that every slice knot is
$(h)$--solvable for all nonnegative half-integers $h$ \cite[Remark
1.3.1]{cot03}. They detect $(n.5)$--solvable knots using the von
Neumann $\rho$--invariants \cite[Theorem 4.2]{cot03}.

Although the von Neumann $\rho$--invariants are defined in a more
general setting, for our purpose, henceforth we study the von
Neumann $\rho$--invariants associated with the representations to
the metabelian group $\Gamma:=\bbq(t)/\Lambda\rtimes \bbz$ where
$\bbz$ is generated by $t$ acting on $\bbq(t)/\Lambda$ by
multiplication. In \cite{cot03} the group $\Gamma$ is called the
\emph{first rationally universal group} and denoted by $\Gamma^U_1$.
Since $\Gamma$ is poly-torsion-free-abelian, the group ring
$\bbz\Gamma$ is a right Ore domain and $\bbq\Gamma$ embeds into its
classical right ring of quotients $\calk_\Gamma$ \cite[Proposition
2.5]{cot03}.

Let $K$ be a knot and $\phi: \pi_1(M_K)\to \Gamma$ a homomorphism.
Then one can define the \emph{von Neumann $\rho$--invariant}
$\rho(M_K,\phi)\in\bbr$ associated with $\phi$ which was introduced
by Cheeger and Gromov \cite{chg85}. When $M_K$ bounds an oriented
compact 4--manifold $W$ with a homomorphism $\psi :\pi_1(W)\to
\Gamma$ extending $\phi$ (i.e., $(M_K,\phi)=\partial(W,\psi)$), the
von Neumann $\rho$--invariant is computed as
$$
\rho(M_K,\phi)=\sigma_\Gamma^{(2)}(W,\psi)-\sigma_0(W)
$$
where $\sigma_\Gamma^{(2)}(W,\psi)$ is the $L^{2}$--signature of the
intersection form on $H_2(W;\calu\Gamma)$ and $\sigma_0(W)$ is the
ordinary signature of $W$. Here $\calu\Gamma$ is the algebra of
(unbounded) operators affiliated to the von Neumann algebra
$\caln\Gamma$ of the group $\Gamma$. Often we simply denote
$\rho(M_K,\phi)$ by $\rho(K,\phi)$. We refer the reader to
\cite[Section 5]{cot03} for more discussion of $L^{2}$--signatures.
The following theorem gives an obstruction for a knot being
$(1.5)$--solvable. Note that since $\Gamma^{(2)} = 0$, $\Gamma$ is
$(1)$--solvable.

\begin{thm}\cite[Theorem 4.2]{cot03}\label{thm:vanishing}
In the above setting, if $W$ is an $(1.5)$--solution for
$K$ then $\rho(K,\phi)=0$.

\end{thm}

In this paper we are interested in the following representations to
$\Gamma$. We construct a representation $\phi_x : \pi_1(M_K) \to
\Gamma$ associated with a given $x\in H_1(M;\Lambda)$ via the
Blanchfield linking form in the following manner. Recall that there
is a nonsingular form called the \emph{Blanchfield linking form}
$$
B\ell :H_1(M_K;\Lambda)\times H_1(M_K;\Lambda)\to \bbq(t)/\Lambda.
$$
Let $\mu$ be a meridian of $K$ which normally generates $\pi_1(M_K)$
and $\epsilon : \pi_1(M_K)\to \bbz$ the abelianization sending $\mu$
to $1$.

\begin{defn}\label{defn:representation}
The representation $\phi_x : \pi_1(M_K)\to \Gamma$ is defined to be
$\phi_x(y)=(B\ell(x,y\mu^{-\epsilon(y)}), \epsilon(y))$ for $y\in
\pi_1(M_K)$,
where $y\mu^{-\epsilon(y)}$ denotes its image in $H_1(M_K;\Lambda)$
as an abuse of notation.
\end{defn}

For a $\Lambda$--submodule $P$ of $H_1(M_K;\Lambda)$, define
$$
P^\perp := \{y\in H_1(M_K;\Lambda) \mid B\ell(x,y)=0 \text{ for all
} x\in P\}.
$$
We say a $\Lambda$--submodule $P$ of
$H_1(M_K;\Lambda)$ is \emph{self-annihilating} if $P=P^\perp$.

\begin{thm}\cite[Theorems 3.6 and 4.4]{cot03}\label{thm:extension}
Let $W$ be a $(1)$--solution for $K$ and
$P:=\ker\{H_1(M_K;\Lambda)\to H_1(W;\Lambda)\}$. Then
\begin{enumerate}
\item $P$ is self-annihilating.
\item $\phi_x$ extends to $\pi_1(W)$ if and only if  $x\in P$.
\end{enumerate}
\end{thm}


\section{Polynomial splitting theorem}

The following is a stronger form of Theorem~\ref{thm:main}. For two
knots $K_1$ and $K_2$, note that the Alexander module of $K_1\#K_2$
is isomorphic with the direct sum of the Alexander modules of $K_1$
and $K_2$.

\begin{thm} \label{thm:splitting}
    Let $K_1$ and $K_2$ be knots and
    let $M_1$, $M_2$, and $M$ be zero surgeries on $K_1$, $K_2$,
    and $K_1\# K_2$, respectively.
    Suppose that the Alexander polynomials $\Delta_{K_1}(t)$ and $\Delta_{K_2}(t)$ are coprime.
    If there exists a self-annihilating submodule $P$ of
    $H_1(M;\Lambda)$ with respect to $B\ell$ such that
    $\rho(K_1\# K_2,\phi_x)=0$ for all $x\in P$,
    then there are self-annihilating submodules $P_i$ of
    $H_1(M_i;\Lambda)$, $i=1,2$, such that $P=P_1\oplus P_2$ and
    $\rho(K_i,\phi_{x_i})=0$ for all $x_i\in P_i$, $i=1,2$.
\end{thm}

\begin{proof}
    First, we show that $P=P_1\oplus P_2$ for some
    $\Lambda$--submodules $P_i\subset H_1(M_i;\Lambda)$, $i=1,2$.
    To show this, we follow the proof of \cite[Lemma~3.1]{kim05}.
    Note that
    $H_1(M;\Lambda)= H_1(M_1;\Lambda)\oplus H_1(M_2;\Lambda)$
    and $B\ell=B\ell_1\oplus B\ell_2$ where $B\ell_1$ and $B\ell_2$ denote
    the Blanchfield linking forms of $M_1$ and $M_2$, respectively.
    We write an element $z\in H_1(M;\Lambda)$ as
    $(x,y)\in H_1(M_1;\Lambda)\oplus H_1(M_2;\Lambda)$.
    Let
    $$
    P_1 := \{x\in H_1(M_1;\Lambda) \mid (x,0)\in P\} \mbox{ and }
    P_2 := \{y\in H_1(M_2;\Lambda) \mid (0,y)\in P\}.
    $$
    Clearly $P_1\oplus P_2\subset P$.
    Conversely, we show $P\subset P_1\oplus P_2$.
    For simplicity, denote $\Delta_{K_i}(t)$ by $\Delta_i$ for $i=1,2$.
    Since $\Delta_1$ and $\Delta_2$ are coprime, they are also coprime
    in $\Lambda$ and hence there are $f$ and $g$
    in $\Lambda$ such that $f\Delta_1 + g\Delta_2 =1$ in $\Lambda$.
    Let $z=(x,y)\in P$.
    Since each $\Delta_i$ annihilates $H_1(M_i;\Lambda)$,
    $\Delta_1 x=0$ and $\Delta_2 y=0$. Thus,
    \begin{align*}
        f\Delta_1 z &= (f\Delta_1 x,f\Delta_1 y)
        = (0, f\Delta_1 y), \\
        g\Delta_2 z &= (g\Delta_2 x,g\Delta_2 y)
        = (g\Delta_2 x, 0).
    \end{align*}
    On the other hand,
    \begin{align*}
        x & = 1x= f\Delta_1 x + g\Delta_2 x=g\Delta_2 x, \\
        y & = 1y= f\Delta_1 y + g\Delta_2 y=f\Delta_1 y.
    \end{align*}
    Thus we have
    \begin{align*}
        (0,y) &=(0,f\Delta_1 y) = f\Delta_1 z, \\
        (x,0) &=(g\Delta_2 x,0) = g\Delta_2 z.
    \end{align*}
    Since $P$ is a $\Lambda$--submodule, we conclude that $(x,0)$ and $(0,y)$ are in $P$ and
    hence $x\in P_1$ and $y\in P_2$.
    Then $z=(x,y)\in P_1\oplus P_2$.

    Next, we will show that each $P_i$ is self-annihilating with respect to $B\ell_i$, i.e.,
    $P_i=P_i^\perp$.
    For any $x_1, x_2\in P_1$, we have $(x_1,0), (x_2,0) \in P$
    hence
    \[
    B\ell_1(x_1,x_2)
    =B\ell_1(x_1,x_2) + B\ell_2(0,0)
    =B\ell( (x_1,0),(x_2,0) )=0.
    \]
    Thus, $P_1\subset P_1^\perp$.
    Conversely, let $x\in P_1^\perp$. For any $z\in P$, we can write
    $z=(x',y)$ for some $x'\in P_1, y\in P_2$. Since $x\in P_1^\perp$,
    $B\ell_1(x,x')=0$ and
    \[
    B\ell( (x,0),(x',y)) = B\ell_1(x,x') + B\ell_2(0,y)=0.
    \]
    Thus, $(x,0)\in P^\perp$. Since $P=P^\perp$, $(x,0)\in P$ and hence $x\in P_1$.
    Thus, we see that $P_1=P_1^\perp$.
    Similarly, we can see that $P_2=P_2^\perp$.

    Finally, we will show that $\rho(K_i,\phi_{x_i})=0$ for all
    $x_i\in P_i$, $i=1,2$.
    We construct a cobordism $C$ between the disjoint union
    $M_1\cup M_2$ and $M$. Though this construction is well-known, it
    is an essential step here and is briefly described below. (For
    details, refer to \cite[Section 4]{cot04}.)
    \begin{figure}
        \setlength{\unitlength}{0.6pt}
        \begin{picture}(291,87)
            \put(0,0){\includegraphics[scale=0.6]{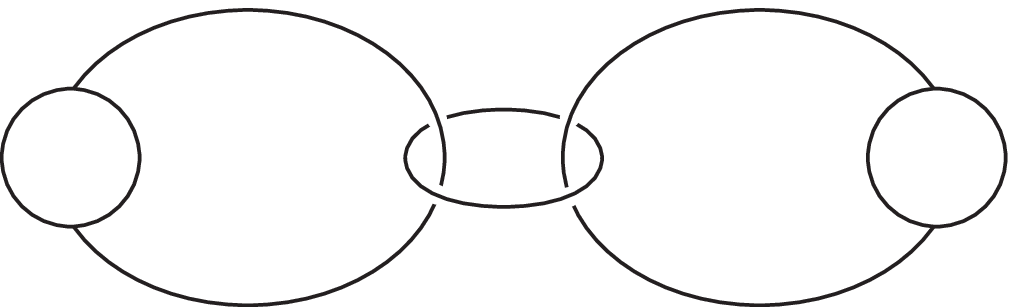}}
            \put(22,43){\makebox(0,0){$K_1$}}
            \put(273,43){\makebox(0,0){$K_2$}}
            \put(145,19){\makebox(0,0){$\eta$}}
            \put(145,62){\makebox(0,0)[b]{\small$0$}}
            \put(73,11){\makebox(0,0){\small$0$}}
            \put(219,11){\makebox(0,0){\small$0$}}
        \end{picture}
        \caption{}
        \label{fig:1}
    \end{figure}
    Attach a 1--handle at the top level
    between $M_1\times [0,1]$ and $M_2\times [0,1]$ so that the
    upper boundary is the zero surgery on the split link $K_1 \cup
    K_2$. Then attach a 2--handle along $\eta$ with zero framing as indicated in Figure~\ref{fig:1}.
    The resulting 4--manifold is the desired $C$. Its boundary at the
    bottom $\partial_{-}C$ is the disjoint union $-(M_1\cup M_2)$.
    To see that $\partial_{+} C = M$, slide
    the 2--handle attached along $K_2$ over the 2--handle attached
    along $K_1$ to get a 2--handle attached along $K_1\# K_2$.
    After the sliding one obtains a surgery diagram where $\eta$ is a meridian of $K_1$ and
    is unlinked from $K_1\#K_2$.
    Using $\eta$ one can unlink and unknot $K_1$ and therefore one can see that
    the zero surgery on $\eta\cup K_1$ is homeomorphic with the
    zero surgery on the Hopf link which is homeomorphic with $S^3$.
    Thus $K_1$ and $\eta$ can be discarded. (Or one can use the
    \emph{slam-dunk} move in \cite[p.163]{gs99}.)
    The result is a surgery diagram for $M$ and this shows that
    $\partial_+ C = M$.

    Let $x\in P_1$. Then $(x,0)\in P$. We obtain the representations
    $\phi_{x}:\pi_1(M_1)\to \Gamma$ and $\phi_{(x,0)}:\pi_1(M)\to \Gamma$
    as defined in
    Definition~\ref{defn:representation}.
    From the above construction, it is easy to see that
    $\pi_1(C)=\pi_1(M_1)\ast \pi_1(M_2)/
    \langle\mu_1\mu_2^{-1}\rangle$, where
    $\ast$, $\langle \ \rangle$, and $\mu_i$ stand for
    the free product, the subgroup normally generated by the given elements,
    and the meridian of $K_i$ at the point that $\eta$ goes over,
    respectively.
    Also, we see that $\pi_1(M)=\pi_1(E_1)\ast\pi_1(E_2)/
    \langle \mu_1\mu_2^{-1}, l_1 l_2 \rangle$, where each $E_i$ is the
    exterior of $K_i$ in $S^3$ and each $l_i$ is the longitude of
    $K_i$.
    Since $\pi_1(C)= \pi_1(M)/\langle l_1\rangle
    = \pi_1(M)/\langle l_2\rangle$ and
    $\phi_{(x,0)}(l_1) = \phi_{(x,0)}(l_2)=0$, the representation
    $\phi_{(x,0)}$ extends over $\pi_1(C)$, i.e.,
    there is a homomorphism $\tilde{\phi}:\pi_1(C)\to\Gamma$ such that
    $\phi_{(x,0)} = \tilde{\phi}\circ i_\ast$, where
    $i_\ast:\pi_1(M) \to \pi_1(C)$ is the homomorphism induced by the
    inclusion $i:M\to C$. From the construction of $C$, it immediately
    follows that $\tilde{\phi}|_{\pi_1(M_1)} = \phi_x$ and
    $\tilde{\phi}|_{\pi_1(M_2)}$ is the zero map.

    By \cite[Lemma 4.2]{cot04} $H_\ast(C;\calk_\Gamma)=0$, and one
    can easily show
    that the inclusion induced map $i_*: H_2(\partial C;\bbz) \to H_2(C;\bbz)$
    is surjective. Therefore
    $\sigma_\Gamma^{(2)}(C,\tilde{\phi})=0 $ and  $\sigma_0(C) = 0$.
    Hence
    \begin{align*}
        \rho(K_1 \# K_2,\phi_{(x,0)}) & =
        \rho(K_1,\tilde{\phi}|_{\pi_1(M_1)}) +
        \rho(K_2,\tilde{\phi}|_{\pi_1(M_2)}) \\
        & = \rho(K_1,\phi_x) + \rho(K_2,0) \\
        & = \rho(K_1,\phi_x).
    \end{align*}
    Here, $\rho(K_2,0)=0$ follows from \cite[Property (2.5)]{cot04}.
    Since $(x,0)\in P$,
    the left hand side is zero by the assumption.
    Therefore $\rho(K_1,\phi_x)=0$ for all $x\in P_1$.
    Similarly, $\rho(K_2,\phi_y)=0$
    for all $y\in P_2$. This completes the proof.
\end{proof}

\section{Examples}
\label{sec:examples}

As an application of Theorem~\ref{thm:splitting}, we present a new
example of an infinite family of knots with vanishing Casson--Gordon
invariants which are linearly independent in $\calc$. In fact, we
will show that they are linearly independent in
$\calf_{1}/\calf_{1.5}$. The authors do not know how to show that
these knots are linearly independent in $\calc$ without using
Theorem~\ref{thm:splitting}.


Let $T$ be the infinite set of positive integers each of which is
divisible by three distinct primes. For $k\in T$, let $\Phi_k(t)$ be
the $k$--th cyclotomic polynomial. As is well-known (for instance
see \cite[Section 2]{tkim04} and \cite[Chapter 5]{cha06}), there is
a slice knot $K'_k$ that has the cyclic rational Alexander module
$H_1(M_{K'_k};\Lambda)=\Lambda/\left( \Phi_k(t)^2 \right)$. Let $J$
be a knot with Arf invariant zero such that the (averaged) integral
of the Levine--Tristram signature function of $J$ is nonzero. For
example, one can take $J$ to be the connected sum of two copies of
the trefoil. Denote this integral by $\rho(J)$. Let $\eta_k$ be the
unknot in the complement of a Seifert surface for $K_k'$ in $S^3$
which represents the homology class generating the rational
Alexander module for $K'_k$ over $\Lambda$, i.e.,
$H_1(M_{K'_k};\Lambda)=\left([\eta_k]\right)$. Such $\eta_k$ exists
since any element in the Alexander module for a knot can be
represented by a simple closed curve in the complement of the knot
which represents a commutator in the knot group and we may assume
that the simple closed curve is unknotted by crossing change (cf.
\cite{rol76}).

We construct a knot $K_k :=K'_k(\eta_k,J)$ using the \emph{satellite
construction} (or \emph{genetic modification} following the
terminology in \cite{cot04}) as follows. Take the union of the
exterior of $\eta_k$ in $S^3$ and the exterior of $J$ in $S^3$ along
their boundary (which is a torus) such that the meridian
(respectively the longitude) of $\eta_k$ is identified with the
longitude (respectively the meridian) of $J$. The resulting ambient
manifold is $S^3$ and $K_k$ is defined to be the image of $K_k'$
under this identification.
In this case, we say that $K_k$ is the \emph{satellite} of the
\emph{companion} $J$ with
the \emph{axis} $\eta_k$ and the \emph{pattern} $K_k'$.
See \cite[Section~3]{cot04} for more details on this construction.

Since $\eta_k$ lies in the complement of a Seifert surface for
$K_k'$, the knots $K_k'$ and $K_k$ have the isomorphic rational
Alexander module. In particular each $K_k$ has the unique nontrivial
proper submodule $ (\Phi_k(t))$. Using this property, in
\cite[Section 6]{tkim04} the second author showed the following:
each $K_k$, $k\in T$, is $(1)$--solvable but not $(1.5)$--solvable.

Furthermore, we prove the linear independence of $K_k$, $k\in T$.

\begin{prop}\label{prop:examples}
    The knots $K_k,\ {k\in T}$, are linearly independent in
    $\calf_{1}/\calf_{1.5}$ and hence in the knot concordance group.
    Moreover, the concordance invariants by Casson--Gordon
    \cite{cg86}, Letsche \cite{let00} and Friedl \cite{fri04} vanish
    on connected sums of copies of $K_k$, $k\in T$.
\end{prop}
\begin{proof}
    As mentioned above, all $K_k$ are $(1)$--solvable. Let $k_1,\ldots,k_l$
    be distinct elements in $T$.
  Suppose that $K:=a_1 K_{k_1}\# \cdots \# a_l K_{k_l}$ (all $a_i$ are integers) is
    $(1.5)$--solvable. We only need to show that $a_i = 0$ for all
    $i$. Suppose to the contrary that $a_i \ne 0$ for some $i$. We may
    assume that
    $a_1 > 0$ by replacing $K_{k_1}$ by $-K_{k_1}$ if necessary. For convenience
    let us denote $K_{k_1}$ by $K_1$.

    Then there is a self-annihilating submodule $P$ of
    $H_1(M_K;\Lambda)$ such that $\rho(M_K,\phi_z)=0$ for all $z\in P$ by
    Theorems~\ref{thm:extension}
    and \ref{thm:vanishing}.
    Since $\Phi_k(t),\ k\in T$, are all irreducible and are pairwise
    coprime, by Theorem~\ref{thm:splitting}
    there exists a self-annihilating submodule $P_1$ of $H_1(M_{a_1 K_1};\Lambda)$ such that
    $\rho(M_{a_1 K_1},\phi_{x})=0$ for all $x\in P_1$.

    Recall that $H_1(M_{a_1 K_1};\Lambda) = \bigoplus\limits^{a_1} H_1(M_{K_1};\Lambda)$.
    Pick a nonzero element $x=(x_1,\ldots,x_{a_1})\in P_1$.
    Suppose that $x_i\in (\Phi_{k_1}(t))$ for all $i=1,\ldots,a_1$. We
    may assume that $x_1 \ne 0$.
    For the Blanchfield linking form
    $B\ell$ of $K_1$,
    $B\ell(\eta_{k_1}, x_i) \ne 0$ if and only if $x_i \ne 0$, since
    $\eta_{k_1}$ generates $H_1(M_{K_1};\Lambda)$ and $B\ell$ is
    nonsingular.
    Define $\epsilon_i$ to be 1 if $x_i\ne 0$ and 0 otherwise, for
    $1\le i \le a_1$. Then
    \begin{align*}
        \rho(M_{a_1 K_1},\phi_x) &
        = \sum_{i=1}^{a_1} \rho(M_{K_{1}},\phi_{x_i}) \\
        & = \sum_{i=1}^{a_1} \epsilon_i\rho(J) \\
        & = \rho(J) + \sum_{i=2}^{a_1} \epsilon_i\rho(J)\\
        & \ne 0,
    \end{align*}
    which contradicts Theorem~\ref{thm:splitting}. The first
    equality can be shown by constructing a cobordism between the
    disjoint union of $a_1$ copies of $M_{K_1}$ and $M_{a_1 K_1}$ as
    we did in the proof of Theorem~\ref{thm:splitting}. Here we
    need to add $(a_1-1)$ 1--handles and the same number of 2--handles
    to construct the cobordism. The argument is almost the same as
    that in the proof of Theorem~\ref{thm:splitting} and hence we do
    not repeat it here. The second equality of the above equation follows from the proof
    of \cite[Proposition 6.4]{tkim04} and Properties $(2.3)$, $(2.4)$ and $(2.5)$ in \cite{cot04}.
    The third equality follows
    since $x_1 \ne 0$ and $\epsilon_1 = 1$.

    Next, if $x_{i_0}\not\in (\Phi_{k_1}(t))$ for some ${i_0}$, $1\le i_0 \le a_1$,
    then since $P_1$ is a $\Lambda$--submodule and $\Phi_{k_1}(t)
    x_{i_0} \ne 0$ in $H_1(M_{K_1};\Lambda)$,
    $$
    \Phi_{k_1}(t) x= (\Phi_{k_1}(t) x_1,\ldots,\Phi_{k_1}(t) x_{a_1})
    $$
    is a nonzero element in
    $P_1$ satisfying $\Phi_{k_1}(t) x_i \in (\Phi_{k_1}(t))$ for all $i$.
    Now, we are back to the previous case leading to a contradiction.

    The second statement follows since
all prime power cyclic branched covers of each $K_k$ are homology
spheres by \cite[Theorem 1.2]{liv02}.
\end{proof}

There are infinitely many knots with vanishing Arf invariant, say
$J_i$, $i\in \bbn$, such that $\rho(J_i)$, $i\in \bbn$, are linearly
independent over $\bbz$ \cite[Proposition 2.6]{cot04}. Using this
the second author \cite{tkim04} showed that for a fixed $k\in T$ the
knots $K^i_k:=K'_k(\eta_k,J_i)$, $i\in \bbn$, are linearly
independent in $\calf_{1}/\calf_{1.5}$. Note that the knots in
Proposition~\ref{prop:examples} cannot be shown to be linearly
independent in this way since they are the satellites of the same
companion $J$. This demonstrates that Theorem~\ref{thm:splitting} is
definitely required to prove Proposition~\ref{prop:examples}.

One can easily prove the following proposition
using the ideas in the proof of Proposition~\ref{prop:examples}, hence we
omit the proof. Note that the family of knots in the following
proposition includes our examples of knots in
Proposition~\ref{prop:examples}.

\begin{prop} The knots $K^i_k$, $i\in \bbn$, $k\in T$, are linearly
independent in $\calf_{1}/\calf_{1.5}$ and hence in $\calc$.
\end{prop}

{\bf Acknowledgements.} The authors thank an anonymous referee for
careful reading of this paper and helpful comments.

\end{document}